\documentclass[a4paper,12pt,twoside]{amsart}

\usepackage[cm]{fullpage}

\usepackage{mymacros,oitp}
\usepackage[pagebackref]{hyperref}

\title{On the Fukaya categories
of projective hypersurfaces of general type}
\author[K.~Ueda]{Kazushi Ueda}
\address{
Graduate School of Mathematical Sciences,
The University of Tokyo,
3-8-1 Komaba,
Meguro-ku,
Tokyo,
153-8914,
Japan.}
\email{kazushi@ms.u-tokyo.ac.jp}
\date{}

\begin{document}

\begin{abstract}
We prove homological mirror symmetry
for projective hypersurfaces of sufficiently high degree
using a functor
from the wrapped Fukaya category of an affine hypersurface
to the Fukaya category of its boundary at infinity.
\end{abstract}

\maketitle

\section{Functors from Lagrangian correspondences}
 \label{sc:intro_functor}

Fix an algebraically closed field $\bfk$ of characteristic zero
as a ground field.
Let $X$ be a smooth hypersurface in $\bC P^{n+1}$
of degree
\begin{align} \label{eq:high degree}
a > 2 n + 1
\end{align}
equipped with the restriction $\omega_X$
of the Fubini--Study form.
% The assumption \pref{eq:high degree}
% implies that
% $X$ is negatively monotone.
Choose
a connection 1-form
$\vartheta$
on the sphere bundle
$
\pi \colon \sfS \det T X \to X
$
associated with the determinant line bundle
of the tangent bundle of $X$
satisfying
$
d \vartheta = - \pi^* \omega_X
$.
An oriented Lagrangian submanifold
$
L
$
of $X$
induces a map
$
\iotat_L \colon L \to \sfS \det T X
$
sending
$
x \in L
$
to (the class of)
$
\det T_x L
$.
An oriented Lagrangian submanifold $L$ is said to be \emph{balanced}
if
$
\iotat_L^* \vartheta
$
is exact
\cite[Section 6]{MR2819674}.
The
$\bZ/2\bZ$-graded
$\bfk$-linear
idempotent-complete
pretriangulated
Fukaya category of $X$
generated by balanced spin Lagrangian submanifolds
is denoted by
$
\fuk{X}
$.

Set
\begin{align} \label{eq:w}
\sfw
&\coloneqq x_1^a + \cdots + x_{n+2}^a,
\end{align}
and let
$
\emf{\bA^{n+2}}{G}{\sfw - x_1 \cdots x_{n+2}}
$
be the idempotent-completion
of the dg category of matrix factorizations
equivariant
under the diagonal action of
\begin{align}
G
\coloneqq
\lc (t_1,\ldots,t_{n+2}) \in (\Gm)^{n+2}
\relmid
t_1^a = \cdots = t_{n+2}^a = t_1 \cdots t_{n+2} = 1 \rc
\end{align}
on $\bA^{n+2}$.

\begin{theorem} \label{th:main}
Under the assumption \pref{eq:high degree},
one has an equivalence
\begin{align}
\fuk{X}
\simeq
\emf{\bA^{n+2}}{G}{\sfw - x_1 \cdots x_{n+2}}.
\end{align}
\end{theorem}

Let
$
Y \subset \bC P^{n+2}
$
be a smooth hypersurface of degree $a$
such that
$
X \subset Y
$
is a smooth hyperplane section.
The complement
$
U \coloneqq Y \setminus X
$
is affine
and hence has a Liouville structure,
whose
% $\bZ/2\bZ$-graded $\bfk$-linear pretriangulated idempotent-complete
wrapped Fukaya category
is denoted by
$\wfuk{U}$.
Let $C$ be the contact boundary
of the Liouville domain
having $U$ as the Liouville completion
such that the Reeb vector field induces an $S^1$-action on $C$
with quotient $X$.
The product
$
C \to U \times X
$
of the inclusion
$
C \hookrightarrow U
$
and the projection
$
C \to X
$
gives a Lagrangian correspondence.
Assume
$n$ is greater than one,
so that both $X$ and $U$ are simply connected.
We will not discuss the case $n=1$,
where homological mirror symmetry is known
by \cite{MR2819674,MR2914956}.
The Lagrangian correspondence $C$ is automatically balanced
since $H^1(C;\bR)$ is zero.
The obstruction to the definition of Floer cohomology of $C$
is given by
\begin{align}
\m_0(1) = \sum_{\beta \in \pi_2(X \times U, C)}
\lb \ev_\beta \rb_!^{\virt}(1)
\end{align}
where $\lb \ev_\beta \rb_!^{\virt}$
is the push-forward
along the evaluation map
\begin{align}
\ev_\beta \colon \cM_1(X \times U, C; \beta) \to C
\end{align}
from the moduli space of pseudoholomorphic disks
with one marked point on the boundary
against the virtual fundamental chain.
The assumption \pref{eq:high degree}
implies
\begin{align}
\deg \lb \ev_\beta \rb_!^{\virt}(1)
&= 2 - \mu(\beta)
\ge 2 + 2 (a - n - 2)
\ge 2 n + 2
> \dim C
\end{align}
and hence $\m_0(1) = 0$,
where $\mu(\beta)$ is the Maslov index of
$
\beta
$.

The unobstructed Lagrangian correspondence $C$
induces an $A_\infty$-functor
\begin{align} \label{eq:Phi}
\Phi \colon \wfuk{U} \to \fuk{X}
\end{align}
by \cite{1706.02131},
which builds on earlier works
including \cite{MR2657646,MR3816496,MR3019714}.
Although \cite{1706.02131} is concerned with compact symplectic manifolds,
the generalization to the Liouville setting is straightforward,
since $C$ is compact.
A less straightforward generalization
to functors associated with non-compact Lagrangian correspondences
is discussed in
\cite{1703.04032,1712.00225}.
A functor
defined by the Lagrangian correspondence $C$
in the monotone case
is also studied in \cite{2406.08852}.

Since compact Lagrangian submanifolds of $U$ are disjoint from the contact boundary,
the functor \pref{eq:Phi}
sends the full subcategory
$\fuk U \subset \wfuk U$
consisting of compact Lagrangian submanifolds
to zero,
and hence descends to a functor
\begin{align} \label{eq:Psi}
\Psi \colon \sfuk U \to \fuk X
\end{align}
from the \emph{stable Fukaya category}
\begin{align}
  \sfuk U \coloneqq \wfuk U / \fuk U.
\end{align}

The Liouville manifold $U$ can be identified
with the Milnor fiber
$
\sfw^{-1}(1)
$
of the polynomial \eqref{eq:w}
defining a Brieskorn--Pham singularity.
It has a $\bZ$-grading
defined by (the tensor square of) the holomorphic volume form
\begin{align}
\Omega \coloneqq \Res \frac{d x_1 \wedge \cdots \wedge d x_{n+2}}{\sfw}.
\end{align}
The $\bZ$-grading of $U$ is unique
since $U$ is simply-connected.
We write the resulting $\bZ$-graded symplectic manifold as $\Ut$
and its $\bZ$-graded wrapped Fukaya category as $\wfuk{\Ut}$.
Let
\begin{align}
K \coloneqq \lc (t_0,\ldots,t_{n+2}) \in (\Gm)^{n+3}
\relmid
t_1^a = \cdots = t_{n+2}^a = t_0 \cdots t_{n+2}
\rc
\end{align}
act diagonally
on
$
\bA^{n+3}
=
\Spec \bfk [x_0,\ldots,x_{n+2}]
$.
% by
% \begin{align}
% (x_0,\ldots,x_{n+2}) \mapsto (t_0 x_0,\ldots,t_{n+2} x_{n+2}).
% \end{align}
By
% \cite[Conjecture 1.4]{MR4442683}
% proved in
\cite[Theorem 4.1]{2406.15915},
one has a quasi-equivalence
\begin{align} \label{eq:Z-graded hms for U}
\wfuk{\Ut}
\simeq
\emf{\bA^{n+3}}{K}{\sfw - x_0 \cdots x_{n+2}}.
\end{align}
A $\bZ/2\bZ$-graded variant of
\pref{eq:Z-graded hms for U}
(cf.~also \cite{MR4713718})
is given by
\begin{align} \label{eq:hms for U}
\wfuk{U}
\simeq
\emf{\bA^{n+3}}{H}{\sfw - x_0 \cdots x_{n+2}},
\end{align}
where
\begin{align}
H \coloneqq \lc (t_0,\ldots,t_{n+2}) \in (\Gm)^{n+3}
\relmid
t_1^a = \cdots = t_{n+2}^a = t_0 \cdots t_{n+2} = 1
\rc.
\end{align}

The contact boundary
% \begin{align}
% \lc (x_1,\ldots,x_{n+2}^a) \in \sfw^{-1}(1) \relmid |x_1|^2 + \cdots + |x_{n+2}|^2 = R \rc
% \end{align}
of the Liouville domain
associated with the Milnor fiber
$
U = \sfw^{-1}(1)
$
is contactomorphic to
\begin{align}
\lc (x_1,\ldots,x_{n+2}) \in \sfw^{-1}(0)
\relmid
|x_1|^2 + \cdots + |x_{n+2}|^2 = 1 \rc,
\end{align}
which has a contact form whose Reeb flow induces the $S^1$-action
\begin{align}
S^1 \ni \zeta
\colon
(x_1,\ldots,x_{n+2}) \mapsto (\zeta x_1,\ldots,\zeta x_{n+2}).
\end{align}
The unit
$
\eta \colon \id_{\wfuk{U}} \to \cap \cup
$
of the adjunction
\begin{align}
(\cup \dashv \cap)
\colon
\wfuk{U}
\overset{\cup}{\underset{\cap}{\rightleftarrows}}
\wfuk{\bC^{n+2}, U}
\end{align}
of the cup functor and the cap functor
fits into the distinguished triangle
\begin{align} \label{eq:Abouzaid--Ganatra}
\mu_* \xto{\ s\ } \id \xto{\ \eta\ } \cap \cup \to \mu_*[1]
\end{align}
by Abouzaid--Ganatra
(cf.~e.g.~\cite{MR4377932}),
where
$
\mu_* \colon \wfuk{U} \to \wfuk{U}
$
is the functor induced by
the clockwise monodromy
$
\mu \colon U \to U
$
of the one-parameter family
$
\bfw \colon \bC^{n+2} \to \bC
$
around the origin,
which is Hamiltonian isotopic to the product of inverse Dehn twists
along the vanishing cycles,
and
$s$
is
the natural transformation
originally
introduced by Seidel \cite{MR2483942}.
Let
$
\chi_0
$
be the character of $K$
sending $(t_0,\ldots,t_{n+2}) \in K$ to $t_0$.
The autoequivalence of
$
\emf{\bA^{n+3}}{K}{\sfw - x_0 \cdots x_{n+2}}
$
induced by the graded variant $\mut_*$ of $\mu_*$
under the equivalence \pref{eq:Z-graded hms for U}
is given by the tensor product
with the dual $\chi_0^\dual$
of $\chi_0$
on the mirrors of vanishing cycles,
and hence on the whole of $\wfuk{\Ut}$
because of \cite[Theorem 6.11]{MR4442683}.
Since the monodromy
$
\mu
$
is trivial outside of a compact set,
the functors
induced
on $\sfuk{U}$
by the functors
$\mu_*$
and
$
\cap \cup
$
are the identity functor
and the zero functor
respectively,
and
the natural transformation $s$
induces an automorphism
of the identity functor of $\sfuk{U}$.
The mirror of the graded variant $\st$ of $s$
gives a natural transformation
from $\chi_0^\dual$ to the identity functor.
The space of natural transformations
from $\chi_0^\dual$ to the identity functor
is the degree $\chi_0$-component
of the extended Hochschild cohomology
in the sense of \cite[Introduction]{MR3270588},
and can be computed to be spanned by $x_0$
as in \cite[Section 3.6]{MR4442683}.
The $a$-th power
$\mut_*^a$
is isomorphic
to the shift functor $[-2(a-n-2)]$,
and $\st^a$ gives an element of
$\HH^{2(a-n-2)} \lb \wfuk{\Ut} \rb$.
The discussion in \cite[Section 5.1]{2406.15915}
shows that
% the stable Fukaya category
$\sfuk{U}$
% which is equivalent to the Rabinowitz Fukaya category
% by \cite{2212.14863,2309.17062},
is equivalent to
$
\wfuk{U} \otimes_{\bfk[s^a]} \bfk[s^{\pm a}]
$,
which is equivalent to
\[
\emf{\bA^{n+3}}{H}{\sfw - x_0 x_1 \cdots x_{n+2}}
\otimes_{\bfk[x_0^a]} \bfk[x_0^{\pm a}]
\simeq
\emf{\bA^{n+2}}{G}{\sfw - x_1 \cdots x_{n+2}}
\otimes_\bfk \bfk[x_0^{\pm a}]
\]
on the mirror side.

The stable Fukaya category
$\sfuk U$
is generated by a set
$
\left\{ L_\bsi \right\}_{\bsi \in \{ 1, \ldots, a-1 \}^{n+1}}
$
of Lefschetz thimbles
along a path to the real infinity on the $x_{n+2}$-axis,
which goes
by $\Psi$
to a distinguished basis
$
\left\{ V_\bsi \right\}_{\bsi \in \{ 1, \ldots, a-1 \}^{n+1}}
$
of vanishing cycles in $X$
(equipped with the zero bounding cochain
% since a bounding cochain must be of degree one
% with respect to the grading such that the Novikov variable
% has degree $2(a-n-2)$, which cannot be non-zero by
because of
\pref{eq:high degree}).
The corresponding object
of
$
\emf{\bA^{n+2}}{G}{\sfw - x_1 \cdots x_{n+2}}
\otimes_\bfk \bfk[x_0^{\pm a}]
$
is
(the stabilization of)
the structure sheaf of the origin of $\bA^{n+2}$
tensored with a character of $G$
and the free $\bfk[x_0^{\pm a}]$-module
of rank one.

Since the zero locus of
$
\sfwv \coloneqq \sfw - x_1 \cdots x_{n+2}
$
has an isolated singularity at the origin,
one has
\begin{align}
\emf{\bA^{n+2}}{G}{\sfwv}
\simeq
\emf{\Spec S}{G}{\sfwv}
\end{align}
by \cite{MR2735755}
where
\begin{align}
S = \bfk \db[ x_1, \ldots, x_{n+2} \db].
\end{align}
Computation in
$
\emf{\Spec S}{G}{\sfwv}
$
using \cite[Theorem 5.9]{MR2824483}
shows that
the
space
% $\bfk[s^{\pm a}]$-module
$
\Hom_{\sfuk U}^* \left( L_\bsi, L_\bsj \right)
$
is non-zero if and only if either
\begin{itemize}
\item
$\bsi \ne \bsj$ and
$
0 \le j_l - i_l \le 1
$
for all
$
l \in \{ 1, \ldots, n+1 \}
$,
where it is a free $\bfk[s^{\pm a}]$-module of rank one
whose generator on the mirror is written as
$
\prod_{l \colon j_l - i_l = 1} \partialbar_l
$
in the notation of \cite{MR2824483},
\item
$\bsi \ne \bsj$ and
$
0 \le i_l - j_l \le 1
$
for all
$
l \in \{ 1, \ldots, n+1 \}
$,
where it is a free $\bfk[s^{\pm a}]$-module of rank one
whose generator on the mirror is
$
\lb \prod_{l \colon i_l - j_l \ne 1} \partialbar_l \rb \partialbar_{n+2}
$,
or
\item
$
\bsi=\bsj
$
where it is a free $\bfk[s^{\pm a}]$-module of rank two
generated by $1$ and
$
\prod_{l=1}^{n+2} \partialbar_l
$.
\end{itemize}

The functor $\Phi$ is full
since it is defined using the Yoneda embedding.
For each Lefschetz thimble $L$,
the functor $\Phi$ sends the morphism
$s_L^a \colon \mu_*^a(L) \to L$
to
an endomorphism
$
\Phi(s_L^a) \colon V \to V
$
of the corresponding vanishing cycle $V$.
Note that the Lagrangian correspondence $C$
and hence the functor $\Phi$ admits
a $\bZ / 2(a-n-2)\bZ$-grading.
Since
the Floer cohomology of the Lagrangian sphere $V$
is concentrated in degrees $0$ and $n$,
which are not congruent modulo $2(a-n-2)$ by \pref{eq:high degree},
the endomorphism $\Phi(s_L^a)$ must be proportional to the identity morphism
of $V$.
Set $\Phi(s_L^a) = \lambda \cdot \id_V$,
where $\lambda \in \bfk$ is non-zero by
(the image by $\Phi$ of)
\pref{eq:Abouzaid--Ganatra},
and independent of $L$
since any pair of vanishing cycles can be connected
by a chain of vanishing cycles
with non-trivial Floer cohomologies between them.
It follows that $\Phi$ induces a full and faithful functor
from
$
\sfuk U \otimes_{\bfk[s^a]} \bfk[s^a]/(s^a - \lambda)
\simeq
\emf{\Spec S}{G}{\sfw - x_1 \cdots x_{n+2}}
$
to $\fuk X$.

To show that $\Phi$ is essentially surjective,
it suffices to show
\begin{align} \label{eq:dim HH^* = dim QH^*}
\dim \HH^*
\lb
\emf{\Spec S}{G}{\sfwv}
\rb
=
\dim \QH^*(X)
\end{align}
by \cite[Theorem 3]{1605.07702}
(see also
\cite{MR4198511}
for a closely related result).
% \footnote{
% \cite{MR4198511,1605.07702}
% depend on properties of the Fukaya category
% announced by Abouzaid--Fukaya--Oh--Ohta--Ono.
% A proof in the positively monotone case is given in \cite{MR3578916},
% and it is straightforward to generalize it
% to the balanced Fukaya category
% of a negatively monotone symplectic manifold.}

For $\gamma = (\alpha_i)_{i=1}^{n+2} \in G$,
we set
\begin{align}
I_\gamma &\coloneqq \lc i \in \{ 1, \ldots, n+2 \} \relmid \alpha_i = 1 \rc, \\
S_\gamma &\coloneqq \bfk \db[ x_i \db]_{i \in I}, \\
\sfwv_\gamma &\coloneqq \sfwv |_{x_i = 0 \text { for } i \nin I}, \\
\Jac_{\sfwv_\gamma}
&\coloneqq
S_\gamma \left/ \lb \frac{\partial \sfwv_\gamma}{\partial x_i} \rb_{i \in I} \right., \\
n_\gamma &\coloneqq \sum_{i \nin I} \chi_i.
\end{align}
It follows from
\cite{MR2824483,MR3084707,MR3108698,MR3270588}
that
\begin{align} \label{eq:HHmf}
\HH^t \lb
\emf{\Spec S}{G}{\sfwv}
\rb
\cong
\bigoplus_{\substack{\gamma \in \ker \chi,\\ t - \dim N_\gamma = 2u }}
\left(
\Jac_{\sfwv_\gamma}
\right)_{u \chi - n_\gamma},
\end{align}
where
$
\left(
\Jac_{\sfwv_\gamma}
\right)_{u \chi - n_\gamma}
$
is the subspace of the Jacobian ring
consisting of elements of degree $u \chi - n_\gamma$.

The summand for $\gamma = \bsone$
% which is known as the untwisted sector
gives
the $\ker \chi$-invariant part
\begin{align}
\lb \bfk \db[ x_1,\ldots,x_{n+2} \db] / (a x_i^{a-1} - x_1 \cdots \hat{x_i} \cdots x_{n+2})_{i=1}^{n+2} \rb^{\ker \chi}
=
\bigoplus_{i=0}^{n} \bfk \ld \lb x_1 \cdots x_{n+2} \rb^i \rd
\end{align}
of the Jacobi ring.
% where
% \begin{align}
% \deg \ld \lb x_1 \cdots x_{n+1} \rb^i \rd = 2 i.
% \end{align}
Note that
\begin{align}
(x_1 \cdots x_{n+2})^{a-1} - (x_1 \cdots x_{n+2})^{n+1}
=
((x_1 \cdots x_{n+2})^{a-n-2}-1) (x_1 \cdots x_{n+2})^{n+1}
\in
\lb \frac{\partial \sfwv}{\partial x_i} \rb_{i=1}^{n+2},
\end{align}
so that
\begin{align}
(x_1 \cdots x_{n+1})^{n+1} = 0 \in \Jac_{\sfwv}.
\end{align}

The summand for
$
\gamma = ( \exp(2 i_k \pi \sqrt{-1}/a) )_{k=1}^{n+2} \ne \bsone
$
is non-zero
if and only if
\begin{align}
1 \le i_k \le a-1
\text{ for all }
k \in \{ 1, \ldots, n+2 \}
\end{align}
and
\begin{align}
i_1 + \cdots + i_{n+2} \equiv 0 \mod a,
\end{align}
where the summand is one-dimensional.
It is known by Griffiths
that $(n-p,p)$-th primitive cohomology of $X$
can be identified with the degree
$((p+1)a-(n+2))$-part
of
the Jacobi ring of
a general homogeneous polynomial $f$ of degree $a$.
If one takes the Fermat polynomial as $f$,
then the Jacobi ring is spanned by
$
 x^\bsi \coloneqq x_1^{i_1} \cdots x_{n+2}^{i_{n+2}}
$
where
$
 0 \le i_k \le a-2
$
for all $k \in \{ 1, \ldots, n+2 \}$.
The set of
$
\bsi \in \{ 0, \ldots, a-2 \}^{n+2}
$
satisfying
\begin{align}
|\bsi| \coloneqq i_1 + \cdots + i_{n+2}
= (p+1)a-(n+2)
\end{align}
is the same as the set of
$
\bsj \coloneqq \bsi + \bsone \in \{ 1, \ldots, a-1 \}^{n+2}
$
satisfying
\begin{align}
|\bsj| = (p+1) a,
\end{align}
and \pref{eq:dim HH^* = dim QH^*} is proved.

\section*{Acknowledgments}

We thank Yank{\i} Lekili
for collaboration on
a series of joint works and
numerous insightful discussions
%collaborations,
% including \cite{2406.15915} in particular,
which led to this work.

\bibliographystyle{amsalpha}
\bibliography{bibs}

\end{document}